\documentstyle[11pt]{article} 

\newtheorem{sect}{}[section]

\newcommand{\vlk}{\mbox{ $v \ell k$}} 
\newread\epsffilein    
\newif\ifepsffileok    
\newif\ifepsfbbfound   
\newif\ifepsfverbose   
\newif\ifepsfdraft     
\newdimen\epsfxsize    
\newdimen\epsfysize    
\newdimen\epsftsize    
\newdimen\epsfrsize    
\newdimen\epsftmp      
\newdimen\pspoints     
\def\epsfbox#1{\global\def\epsfllx{72}\global\def\epsflly{72}%
   \global\def\epsfurx{540}\global\def\epsfury{720}%
   \def\lbracket{[}\def\testit{#1}\ifx\testit\lbracket
   \let\next=\epsfgetlitbb\else\let\next=\epsfnormal\fi\next{#1}}%
\def\epsfgetlitbb#1#2 #3 #4 #5]#6{\epsfgrab #2 #3 #4 #5 .\\%
   \epsfsetgraph{#6}}%
\def\epsfnormal#1{\epsfgetbb{#1}\epsfsetgraph{#1}}%
\def\epsfgetbb#1{%
\openin\epsffilein=#1
\ifeof\epsffilein\errmessage{I couldn't open #1, will ignore it}\else
   {\epsffileoktrue \chardef\other=12
    \def\do##1{\catcode`##1=\other}\dospecials \catcode`\ =10
    \loop
       \read\epsffilein to \epsffileline
       \ifeof\epsffilein\epsffileokfalse\else
          \expandafter\epsfaux\epsffileline:. \\%
       \fi
   \ifepsffileok\repeat
   \ifepsfbbfound\else
    \ifepsfverbose\message{No bounding box comment in #1; using defaults}\fi\fi
   }\closein\epsffilein\fi}%
\def\epsfclipoff{\def\epsfclipstring{\ifepsfdraft\space clip\fi}}%
\epsfclipoff
\def\epsfsetgraph#1{%
   \epsfrsize=\epsfury\pspoints
   \advance\epsfrsize by-\epsflly\pspoints
   \epsftsize=\epsfurx\pspoints
   \advance\epsftsize by-\epsfllx\pspoints
   \epsfxsize\epsfsize\epsftsize\epsfrsize
   \ifnum\epsfxsize=0 \ifnum\epsfysize=0
      \epsfxsize=\epsftsize \epsfysize=\epsfrsize
      \epsfrsize=0pt
     \else\epsftmp=\epsftsize \divide\epsftmp\epsfrsize
       \epsfxsize=\epsfysize \multiply\epsfxsize\epsftmp
       \multiply\epsftmp\epsfrsize \advance\epsftsize-\epsftmp
       \epsftmp=\epsfysize
       \loop \advance\epsftsize\epsftsize \divide\epsftmp 2
       \ifnum\epsftmp>0
          \ifnum\epsftsize<\epsfrsize\else
             \advance\epsftsize-\epsfrsize \advance\epsfxsize\epsftmp \fi
       \repeat
       \epsfrsize=0pt
     \fi
   \else \ifnum\epsfysize=0
     \epsftmp=\epsfrsize \divide\epsftmp\epsftsize
     \epsfysize=\epsfxsize \multiply\epsfysize\epsftmp
     \multiply\epsftmp\epsftsize \advance\epsfrsize-\epsftmp
     \epsftmp=\epsfxsize
     \loop \advance\epsfrsize\epsfrsize \divide\epsftmp 2
     \ifnum\epsftmp>0
        \ifnum\epsfrsize<\epsftsize\else
           \advance\epsfrsize-\epsftsize \advance\epsfysize\epsftmp \fi
     \repeat
     \epsfrsize=0pt
    \else
     \epsfrsize=\epsfysize
    \fi
   \fi
   \ifepsfverbose\message{#1: width=\the\epsfxsize, height=\the\epsfysize}\fi
   \epsftmp=10\epsfxsize \divide\epsftmp\pspoints
   \vbox to\epsfysize{\vfil\hbox to\epsfxsize{%
      \ifnum\epsfrsize=0\relax
        \includegraphics{\ifepsfdraft}%
      \else
        \epsfrsize=10\epsfysize \divide\epsfrsize\pspoints
        \includegraphics{\ifepsfdraft}%
      \fi
      \hfil}}%
\global\epsfxsize=0pt\global\epsfysize=0pt}%
{\catcode`\%=12 \global\let\epsfpercent=
\long\def\epsfaux#1#2:#3\\{\ifx#1\epsfpercent
   \def\testit{#2}\ifx\testit\epsfbblit
      \epsfgrab #3 . . . \\%
      \epsffileokfalse
      \global\epsfbbfoundtrue
   \fi\else\ifx#1\par\else\epsffileokfalse\fi\fi}%
\def\epsfempty{}%
\def\epsfgrab #1 #2 #3 #4 #5\\{%
\global\def\epsfllx{#1}\ifx\epsfllx\epsfempty
      \epsfgrab #2 #3 #4 #5 .\\\else
   \global\def\epsflly{#2}%
   \global\def\epsfurx{#3}\global\def\epsfury{#4}\fi}%
\def\epsfsize#1#2{\epsfxsize}
\begin{document}

\title{Quandle Homology Groups, Their Betti Numbers, and Virtual Knots}

\author{
J. Scott Carter \\
University of South Alabama \\
Mobile, AL 36688 \\ carter@mathstat.usouthal.edu \and
Daniel Jelsovsky \\
University of South Florida \\
Tampa, FL 33620 \\ jelsovsk@math.usf.edu \and
Seiichi Kamada \\
Osaka City University \\
Osaka 558-8585, JAPAN\\ kamada@sci.osaka-cu.ac.jp \\
skamada@mathstat.usouthal.edu
\and
Masahico Saito \\
University of South Florida \\
Tampa, FL 33620 \\ saito@math.usf.edu
}
\maketitle
\begin{abstract}
Lower bounds of betti numbers for  homology  groups 
of racks and quandles will be given using the quotient homomorphism
to the orbit quandles. Exact sequences relating various types of homology 
groups are  analyzed. Geometric methods of proving non-triviality
of cohomology groups are also given, using virtual knots. 
The results can be applied to knot theory as the first step 
towards evaluationg the state-sum invariants defined from
quandle cohomology.
\end{abstract}

\section{Introduction}

In \cite{CJKLS}, the authors and L.~Langford
introduced
a notion of cohomology groups of a quandle  to define
a state-sum invariant (the CJKLS invariant) 
of knotted curves and knotted surfaces. A  similar notion
for racks had been defined by R.~Fenn, C.~Rourke and B.~Sanderson \cite{FRS1}.
One of the
purposes of this paper is to relate these two homology theories. To this end, 
we will
define a short
exact sequence of chain complexes associated with a quandle and define
three kinds of
homology (and cohomology) groups of the quandle.
A second purpose is  to give a
lower bound
on the Betti numbers of the three kinds of homology groups. This helps us
to determine
non-triviality of the homology groups of a quandle. The lower bound is
valid also for a
rack if the homology is in the sense of \cite{FRS1}. In this case, 
the methods generalize
an idea of  Greene \cite{Greene} called orbit-writhe.
 A third  purpose is to illustrate
geometric techiques that
use the CJKLS invariants and generalize some of Greene's methods. 
These techniques will also 
demonstrate that large classes of quandles have non-trivial homology.
Since one needs non-trivial cocycles to define the CJKLS invariants,
non-triviality of (co)homology groups provides the first step towards
obtaining the invariants.

We thank Dan Silver for some useful comments.

\section{Basic Notions}

A {\it quandle\/}, $X$, is a set with a binary operation
$\ast$ such that 

\noindent
(I. {\sc idempotency}) for any $a \in X$, $a \ast a=a$,

\noindent
(II. {\sc right-invertibility}) for any $a, b \in X$,
 there is a unique $c \in X$ such that $ a =
c\ast b$,
and 

\noindent 
(III. {\sc self-distributivity}) for any $a, b, c \in X$, we have
 $(a \ast b) \ast c =
(a \ast c) \ast (b \ast c)$, cf.~\cite{Joyce}.
A {\it rack\/} is a set with a binary operation that satisfies (II) and
(III), cf.~\cite{FR}.
A similar notion is known as an {\it automorphic set\/}, cf.~\cite{B}.

\begin{sect}{\bf Examples of quandles.\ }{\rm 
Any set $X$ with the operation $x*y=x$ for any $x,y \in X$ is 
a quandle called the {\it trivial} quandle. 
The trivial quandle of $n$ elements is denoted by $T_n$.

Any group is a quandle by conjugation as operation.
Any subset that is closed under conjugation is also a quandle.
For example, the set, $QS(5)$, of non-identity elements of the permutation 
group on $3$ letters is a quandle. 

Let $n$ be a positive integer. 
For elements  $i, j \in \{ 0, 1, \ldots , n-1 \}$, define
$i\ast j= 2j-i$ where the sum on the right is reduced mod $n$. 
Then $\ast$ defines a quandle structure  called the {\it dihedral quandle},
 $R_n$.
This set can be identified with  the set of reflections of a regular $n$-gon
 with conjugation
as the quandle operation.

Any $\Lambda={\bf Z}[T, T^{-1}]$-module $M$ is a quandle with 
$a*b=Ta+(1-T)b$, $a,b \in M$, called an {\it  Alexander  quandle}. 
Furthermore for a positive integer $n$, a {\it mod-$n$ Alexander  quandle}
${\bf Z}_n[T, T^{-1}]/(h(T))$
is a quandle 
for 
a Laurent polynomial $h(T)$.
The mod-$n$ Alexander quandle is finite 
if the coefficients of the  
highest and lowest degree terms 
of $h$  
 are $\pm 1$.

See \cite{B}, \cite{FR}, \cite{Joyce}, or \cite{Matveev}  for further examples.

}\end{sect}

\begin{sect}{\bf  Homomorphisms and orbits. \ }
\label{hamandeggs}
{\rm
A function $f: X \rightarrow  Y$ between quandles
or racks  is a {\it homomorphism}
if $f(a \ast b) = f(a) * f(b).$ Given a quandle homomorphism, $f$, define
for $x \in X$,
$$E_x=[x]=\{ y \in X | f(x)=f(y) \}= f^{-1} (f(x)).$$
The set $E_x$ is called the  {\it equalizer} of $x$; it is a subquandle of 
$X$. The equalizers form a partition or equivalence relation
$\equiv$ 
 on $X$. Clearly, 
$X / \! \equiv$ is a quandle isomorphic to the image of $f$. If $f$ 
is surjective, then
the quandle $Y$ is said to be a {\it quotient quandle}.

Let $X$ denote a quandle.
{}From Axiom~II, each element $b \in X$ defines a bijection
$S(b) : X \to X$ with $aS(b) = a \ast b$. The bijection is an automorphism
by Axiom~III.
For a word $w = b_1^{\epsilon_1} \dots b_n^{\epsilon_n}$ 
where
$b_1, \dots, b_n
\in X;
\epsilon_1, \dots, \epsilon_n \in \{\pm 1\}$,
we define
$a \ast w = aS(w)$ by
$aS(b_1)^{\epsilon_1}\dots S(b_n)^{\epsilon_n}$.
An automorphism of $X$ is called an {\it inner-automorphism\/} of $X$ if it
is $S(w)$ for a
word $w$. (The notation $S(b)$ follows Joyce's paper \cite{Joyce} and $a
\ast w$ ($= a^w$)
follows Fenn-Rourke \cite{FR}.)

We define a relation on $X$ by $ a \sim b$ if
$a$ is mapped to $b$ by an
inner-automorphism of $X$. The relation $\sim$ is an equivalence relation.
The {\it orbit\/} of $a \in X$ is the
equivalence class of $a$, which is denoted by ${\rm Orb}(a)$.
The set of equivalence classes of $X$ by $\sim$ is denoted by ${\rm Orb}(X)$.
When we regard ${\rm Orb}(X)$ as a trivial quandle,
the projection map $\pi: X \to {\rm Orb}(X)$ is a 
quandle
homomorphism. 
In this case, ${\rm Orb}(X)$ is called the {\it orbit quandle}  of $X$.

For  $a \in X$, the {\it weak orbit} \cite{Joyce} of  $a $  is
$\{ f(a) |  f  $ is an automorphism of  $X \}$.
The 
orbit   of  $a$  is
$\{ f(a) |  f  $ is an inner-automorphism of $X \}$.
A quandle is {\it weakly homogeneous} \cite{Joyce}
if it 
 has only one weak orbit.
A quandle is {\it  homogeneous} 
if is has only one  orbit.
A quandle homomorphism $f:X\rightarrow Y$ is said to be {\it locally-homogeneous} if each equalizer, $E_a$ is a homogenous quandle.
}
\end{sect}

\begin{sect} {\bf Lemma.\/} 
Let $h: X \rightarrow Y$ be a   homomorphism.

$(1)$ If $a,b \in X$ are in the same orbit, then $E_a$ and
$E_b$ are isomorphic.

$(2)$  If $h$ is surjective and $Y$ is homogeneous, then
for any $a,b \in X$
 the subquandles $E_{a}$ and $E_{b}$ are isomorphic.
\end{sect}
{\it Proof.\/} We prove (2); a simliar argument gives (1).
Let $a,b \in X$, let $x=f(a)$ and $y=f(b)$. Since  
$Y$ is homogeneous  there is a word, $w$, in the  free group on
$Y$ such that $x=y\ast w$. Say $w=y_1^{\epsilon 1} \cdots y_n^{\epsilon n}$.
Choose preimages $x_i$ for each of the $y_i$, and define a word 
$v= x_1^{\epsilon 1} \cdots x_n^{\epsilon n}$ in the free group on $X$.
Then the inner automorphism $x \mapsto x*v$ defined on $X$ when restricted to 
$E_a$ is an isomorphism onto $E_b$. 
$\Box$

\section{ Homology and Cohomology}
Let $C_n^{\rm R}(X)$ be the free abelian group generated by
$n$-tuples $(x_1, \dots, x_n)$ of elements of a rack/quandle $X$. Define a
homomorphism
$\partial_{n}: C_{n}^{\rm R}(X) \to C_{n-1}^{\rm R}(X)$ by \begin{eqnarray}
\lefteqn{
\partial_{n}(x_1, x_2, \dots, x_n) } \nonumber \\ && =
\sum_{i=2}^{n} (-1)^{i}\left[ (x_1, x_2, \dots, x_{i-1}, x_{i+1},\dots, x_n) \right.
\nonumber \\
&&
- \left. (x_1 \ast x_i, x_2 \ast x_i, \dots, x_{i-1}\ast x_i, x_{i+1}, \dots, x_n) \right]
\end{eqnarray}
for $n \geq 2$ 
and $\partial_n=0$ for 
$n \leq 1$. 
 Then
$C_\ast^{\rm R}(X)
= \{C_n^{\rm R}(X), \partial_n \}$ is a chain complex.

Let $C_n^{\rm D}(X)$ be the subset of $C_n^{\rm R}(X)$ generated
by $n$-tuples $(x_1, \dots, x_n)$
with $x_{i}=x_{i+1}$ for some $i \in \{1, \dots,n-1\}$ if $n \geq 2$;
otherwise let $C_n^{\rm D}(X)=0$. If $X$ is a quandle, then
$\partial_n(C_n^{\rm D}(X)) \subset C_{n-1}^{\rm D}(X)$ and
$C_\ast^{\rm D}(X) = \{ C_n^{\rm D}(X), \partial_n \}$ is a sub-complex of
$C_\ast^{\rm
R}(X)$. Put $C_n^{\rm Q}(X) = C_n^{\rm R}(X)/ C_n^{\rm D}(X)$ and 
$C_\ast^{\rm Q}(X) = \{ C_n^{\rm Q}(X), \partial'_n \}$,
where $\partial'_n$ is the induced homomorphism.
Henceforth, all boundary maps will be denoted by $\partial_n$.

For an abelian group $G$, define the chain and cochain complexes
\begin{eqnarray}
C_\ast^{\rm W}(X;G) = C_\ast^{\rm W}(X) \otimes G, \quad && \partial =
\partial \otimes {\rm id}; \\ C^\ast_{\rm W}(X;G) = {\rm Hom}(C_\ast^{\rm
W}(X), G), \quad
&& \delta= {\rm Hom}(\partial, {\rm id})
\end{eqnarray}
in the usual way, where ${\rm W} = {\rm R}$ if $X$ is a rack,
or one of ${\rm D}$, ${\rm R}$, ${\rm Q}$ if $X$ is a quandle.

\begin{sect}{\bf Definition.\/} {\rm
The $n$\/th {\it rack homology group\/} and the $n$\/th {\it rack
cohomology group\/} \cite{FRS1}
of a rack/quandle $X$ with coefficient group $G$ are \begin{eqnarray}
H_n^{\rm R}(X;G) = H_{n}(C_\ast^{\rm R}(X;G)), \quad
H^n_{\rm R}(X;G) = H^{n}(C^\ast_{\rm R}(X;G)). \end{eqnarray}
The $n$\/th {\it degeneration homology group\/} and the $n$\/th
 {\it degeneration cohomology group\/}
 of a quandle $X$ with coefficient group $G$ are
\begin{eqnarray}
H_n^{\rm D}(X;G) = H_{n}(C_\ast^{\rm D}(X;G)), \quad
H^n_{\rm D}(X;G) = H^{n}(C^\ast_{\rm D}(X;G)). \end{eqnarray}
The $n$\/th {\it quandle homology group\/}  and the $n$\/th
{\it quandle cohomology group\/ } \cite{CJKLS} of a quandle $X$ with coefficient group $G$ are
\begin{eqnarray}
H_n^{\rm Q}(X;A) = H_{n}(C_\ast^{\rm Q}(X;G)), \quad
H^n_{\rm Q}(X;A) = H^{n}(C^\ast_{\rm Q}(X;G)). \end{eqnarray}

The homology group of a rack in the sense of \cite{FRS1} is $H_n^{\rm R}(X;G)$
and the cohomology of a quandle used in \cite{CJKLS} is $H^n_{\rm Q}(X;A)$.
Refer to
\cite{FRS1}, \cite{FRS2}, \cite{Flower}, \cite{Greene} for some calculations and
applications of the rack homology groups, and to \cite{CJKLS}, \cite{CJKS}
for those of quandle cohomology groups.

\begin{sloppypar}
The cycle and boundary groups 
(resp. cocycle and coboundary groups)
are denoted by $Z_n^{\rm W}(X;G)$ and $B_n^{\rm W}(X;G)$
(resp.  $Z^n_{\rm W}(X;G)$ and $B^n_{\rm W}(X;G)$), 
 so that
$$H_n^{\rm W}(X;G) = Z_n^{\rm W}(X;G)/ B_n^{\rm W}(X;G),
\; H^n_{\rm W}(X;G) = Z^n_{\rm W}(X;G)/ B^n_{\rm W}(X;G)$$
where ${\rm W}$ is one of ${\rm D}$, ${\rm R}$, ${\rm Q}$.
We will omit the coefficient group $G$ if $G = {\bf Z}$ as usual. We denote by
$\beta_n^{\rm W}(X)$ the {\it Betti numbers\/} of $X$ determined by the
homology group
$H_n^{\rm W}(X)$.\end{sloppypar}
}\end{sect}

\begin{sect}{\bf Lemma.\/}
If $X= X_m$ is a finite rack of $m$ elements, then the ranks of the free
abelian groups
$C_n^{\rm D}(X)$, $C_n^{\rm R}(X)$, $C_n^{\rm Q}(X)$ are given by
\begin{eqnarray}
{\rm rank}C_n^{\rm D}(X_m) = a_n, \quad {\rm rank}C_n^{\rm R}(X_m) = m^n, \quad
{\rm rank}C_n^{\rm Q}(X_m) = b_n, \end{eqnarray}
where $a_n+b_n = m^n$ and $b_n = m(m-1)^{n-1}$ for $n\ge 1.$
\end{sect}
 {\it Proof.} We prove that 
\begin{eqnarray}
a_1=0, \quad a_n = (m-1)a_{n-1} + m^{n-1} \quad (n \geq 2), \end{eqnarray}
by induction on $n$.
 By definition,
$C_1^{\rm D}(X)=0$ and $a_1=0$. The number of $n$-tuples $(x_1, \dots,
x_n)$ with
$x_{n-1}=x_n$ is $m^{n-1}$. By induction hypothesis, the number of
$(n-1)$-tuples $(x_1,
\dots, x_{n-1})$ with $x_i=x_{i+1}$ for some $i$ is $a_{n-1}$. So the number of
$n$-tuples $(x_1, \dots, x_n)$ such that $x_i=x_{i+1}$ for some $i$ and
$x_{n-1} \neq
x_n$ is $a_{n-1} \times (m-1)$. Thus we have $a_n= (m-1)a_{n-1} + m^{n-1}$
for $n \geq
2$. By definition we have $b_n = m^n - a_n$ and hence 
\begin{eqnarray}
b_1 =m, \quad b_n= (m-1)b_{n-1} \quad (n \ge 2). 
\end{eqnarray}
Solving this recursion, we have $b_n= m(m-1)^{n-1}$ for $n\ge 1$.
$\Box$

For example,
\begin{eqnarray}
{\rm rank}C_n^{\rm D}(X_2) = 2^n-2, \quad {\rm rank}C_n^{\rm R}(X_2) = 2^n,
\quad
{\rm rank}C_n^{\rm Q}(X_2) = 2.
\end{eqnarray}
Thus, in general, calculation of the quandle homology of a finite quandle is
easier than calculation of the rack homology if one calculates them directly
from the definition.

Let $f : X \to Y$ be a rack homomorphism. It induces a chain map
$f_\sharp : C_\ast^{\rm W}(X) \to C_\ast^{\rm W}(Y)$ in the natural way, and
homomorphisms
$f_\ast : H_n^{\rm W}(X;G) \to H_n^{\rm W}(Y;G)$ and
$f^\ast : H^n_{\rm W}(Y;G) \to H^n_{\rm W}(X;G)$, where ${\rm W} = {\rm R}$ if $X, Y$
are racks, or one of ${\rm D}$, ${\rm R}$, ${\rm Q}$ if $X, Y$ are
quandles. They are
called the {\it homomorphisms induced from $f$.\/}

\begin{sect}\label{sect:BasicHomologyLongExact}
{\bf Proposition (Basic Homology Long Exact Sequence).\/}
If $X$ is a quandle, there is
a long exact sequence \begin{eqnarray}
\cdots \stackrel{\partial_\ast}{\to} H_n^{\rm D}(X;G) \stackrel{i_\ast}{\to}
H_n^{\rm R}(X;G)
\stackrel{j_\ast}{\to} H_n^{\rm Q}(X;G)
\stackrel{\partial_\ast}{\to} H_{n-1}^{\rm D}(X;G) \to \cdots
\end{eqnarray}
which is natural with respect to homomorphisms induced from quandle
homomorphisms.
\end{sect}
{\it Proof.\/}
For each $n$ the following short exact sequence is split. \begin{eqnarray}
0 \to C_n^{\rm D}(X) \stackrel{i}{\to} C_n^{\rm R}(X) \stackrel{j}{\to}
C_n^{\rm Q}(X) \to 0.
\end{eqnarray}
So we have an exact sequence of chain complexes \begin{eqnarray}
0 \to C_\ast^{\rm D}(X)\otimes G \stackrel{i}{\to} C_\ast^{\rm R}(X)\otimes G
\stackrel{j}{\to} C_\ast^{\rm Q}(X)\otimes G \to 0 \end{eqnarray} that
induces the long
exact sequence on homology. $\Box$

\begin{sect}\label{sect:UnivCoefThm}
{\bf Proposition (Universal Coefficient Theorem).\/}
There exist split exact sequences
\begin{eqnarray}
0 \to H_n^{\rm W}(X) \otimes G \to H_n^{\rm W}(X; G) \to
{\rm Tor}(H_{n-1}^{\rm W}(X), G) \to 0 \\
0 \to {\rm Ext}(H_{n-1}^{\rm W}(X), G) \to H^n_{\rm W}(X; G) \to
{\rm Hom}(H_n^{\rm W}(X), G) \to 0,
\end{eqnarray}
where ${\rm W} = {\rm R}$ if $X$ is
a rack, or one of ${\rm D}$, ${\rm R}$, ${\rm Q}$ if $X$ is a quandle.
\end{sect}
{\it Proof.\/} Since $\{C_n^{\rm W}(X)\}$ is a chain complex of free
abelian groups,
we have the result. $\Box$

By the universal coefficient theorem, it is sufficient to know the homology
groups with
integer coefficients. So we will investigate the basic homology long exact
sequence with
$G = {\bf Z}$.

\begin{sect}{\bf Example.\/} 
{\rm
(1) Let $R_3$ be the dihedral quandle of three elements. By a direct
calculation
from the definition, we have
\begin{eqnarray}
H_1^{\rm Q}(R_3) = {\bf Z}, \quad H_2^{\rm Q}(R_3) = 0. \end{eqnarray}
Thus we have that
\begin{eqnarray}
H_2^{\rm Q}(R_3;G) = 0, \quad H^2_{\rm Q}(R_3;G) = 0 \end{eqnarray}
for any coefficient group $G$.

(2) Let $R_4$ be the dihedral quandle of four elements. By a direct calculation
from the definition, we have
\begin{eqnarray}
H_1^{\rm Q}(R_4) = {\bf Z}^2, \quad
H_2^{\rm Q}(R_4) = {\bf Z}^2 \oplus ({\bf Z}_2)^2. \end{eqnarray}
Thus we have that
\begin{eqnarray}
H_2^{\rm Q}(R_4;Z_2) = ({\bf Z}_2)^4, \quad H^2_{\rm Q}(R_4;Z_2) = ({\bf Z}_2)^4
\end{eqnarray}
and
\begin{eqnarray}
H_2^{\rm Q}(R_4;Z_m) = ({\bf Z}_m)^2, \quad H^2_{\rm Q}(R_4;Z_m) = ({\bf Z}_m)^2
\end{eqnarray}
for any positive odd integer $m$.
} \end{sect}

\begin{sect}{\bf Trivial quandle.\/} {\rm
Let $T_m$ be the trivial quandle with $m$ ($ < \infty$) elements. Since
$\partial_n :
C_n^{\rm R}(T_m) \to C_{n-1}^{\rm R}(T_m)$ is the $0$-map, the boundary operators
$\partial_\ast$ in the basic homology long exact sequence (with $G={\bf Z}$) are
$0$-maps and it is decomposed into the short exact sequences
\begin{eqnarray}
0 \to H_n^{\rm D}(T_m) \to H_n^{\rm R}(T_m) \to H_n^{\rm Q}(T_m) \to 0,
\end{eqnarray}
which are identified with the short exact sequences \begin{eqnarray}
0 \to C_n^{\rm D}(T_m) \to C_n^{\rm R}(T_m) \to C_n^{\rm Q}(T_m) \to 0.
\end{eqnarray}
In particular, we have
\begin{eqnarray}
\beta_n^{\rm D}(T_m) = a_n, \quad
\beta_n^{\rm R}(T_m) = m^n, \quad
\beta_n^{\rm Q}(T_m) = b_n
\label{eqn:BettiTrivial}
\end{eqnarray}
where $a_n$ and $b_n$ are as before.
} \end{sect}

For simplicity,
we assume that $|{\rm Orb}(X)| = m < \infty$ in what follows.

Let $\pi: X \to {\rm Orb}(X) = T_m$ be the projection
from a quandle $X$ to its orbit quandle identified with $T_m$.
{}From the naturality of the
basic homology long exact sequence, we have a commutative diagram
\begin{eqnarray} \label{eqn:BasicTrivializingDiagram}
\begin{array}{ccccccccc}
\cdots & \stackrel{\partial_\ast}{\to}
& H_n^{\rm D}(X) & \stackrel{i_\ast}{\to}
& H_n^{\rm R}(X) & \stackrel{j_\ast}{\to}
& H_n^{\rm Q}(X) &
\stackrel{\partial_\ast}{\to} & \cdots
\\ {} & {} & \downarrow & {} & \downarrow & {} & \downarrow & {} & {} \\
0 & \to
& H_n^{\rm D}(T_m) & \stackrel{i_\ast}{\to}
& H_n^{\rm R}(T_m) & \stackrel{j_\ast}{\to}
& H_n^{\rm Q}(T_m) & \to & 0, \end{array}
\end{eqnarray}
where the vertical maps are the induced homomorphisms $\pi_\ast$.

\begin{sect}{\bf Remark (Orbit-Writhe).\/} {\rm
Let $\pi: X \to {\rm Orb}(X) = T_m$ be the projection from a quandle $X$
to its orbit quandle. $C_n^{\rm R}(T_m)$ is freely generated by $n$-tuples
$\vec{\omega}= (\omega_1,
\dots, \omega_n)$ of elements of $T_m = {\rm Orb}(X)$. Let $\vec{\omega}$ be one of
the generators,
and let $p_{\vec{\omega}}: H_n^{\rm R}(T_m)= C_n^{\rm R}(T_m) \to {\bf Z}$ be the
projection to the
factor generated by $\vec{\omega}$. The composition
$p_{\vec{\omega}} \circ \pi_\ast : H_n^{\rm R}(X) \to {\bf Z}$ or $Z_n^{\rm R}(X) \to
H_n^{\rm R}(X) \to {\bf Z}$ is the {\em $\vec{\omega}$-orbit writhe\/} in the sense of
Greene
\cite{Greene}. } \end{sect}

\begin{sect}\label{lem:FirstHomology}
{\bf Proposition. \ }
For a quandle $X$, $H_1^{\rm D}(X)= 0$.
$H_1^{\rm R}(X)$ and $H_1^{\rm Q}(X)$ are free abelian groups of rank
$m =|{\rm Orb}(X)|$.
\end{sect}
{\it Proof.\/}
By definition, $H_1^{\rm D}(X)= H_0^{\rm D}(X)=0$ for any quandle $X$. By
the basic
homology long exact sequence, $H_1^{\rm R}(X)$ is isomorphic to $H_1^{\rm
Q}(X)$. The cycle
group $Z_1^{\rm R}(X)$ is freely generated by elements of $X$, and the
boundary group
$B_1^{\rm R}(X)$  is generated by the images
$\partial_2((x,y))= (x) - (x \ast y)$ for all pairs $(x,y)$ of the elements
of $X$.
Therefore if $x \sim y$, then $[x] = [y]$ in $H_1^{\rm R}(X)$. 
Hence $H_1^{\rm R}(X)$ is
generated by $\{ [x_\omega] | \omega \in {\rm Orb}(X) \}$,
 where $x_\omega$ is a
representative of an orbit $\omega$ in ${\rm Orb}(X)$. In the
diagram~(\ref{eqn:BasicTrivializingDiagram}) 
with $n=1$, $H_1^{\rm R}(T_m)$
is the free
abelian group generated by $\{ [\omega] | \omega \in {\rm Orb}(X)=T_m \}$,
and $\pi_\ast:
H_1^{\rm R}(X) \to H_1^{\rm R}(T_m)$ maps $[x_\omega]$ to $[\omega]$.
Therefore $\pi_\ast:
H_1^{\rm R}(X) \to H_1^{\rm R}(T_m)$ is an isomorphism. $\Box$

\begin{sect}\label{lem:SecondHomology}
{\bf Proposition. \ }
For a quandle $X$, $H_2^{\rm D}(X)$ is a free abelian group
of rank $m =|{\rm Orb}(X)|$. The boundary operator
$\partial_\ast: H_3^{\rm Q}(X) \to
H_2^{\rm D}(X)$ is the $0$-map. Hence the basic homology long exact sequence
has a short
exact factor
\begin{eqnarray}
0 \to H_2^{\rm D}(X) & \stackrel{i_\ast}{\to} H_2^{\rm R}(X)
\stackrel{j_\ast}{\to}
H_2^{\rm Q}(X) \to 0. \end{eqnarray}
\end{sect}
{\it Proof.\/}
$Z_2^{\rm D}(X)= C_2^{\rm D}(X)$, which is generated by $(x,x)$
for all $x \in X$. $B_2^{\rm D}(X)$ is generated by
$\partial_3((x,x,y))= - (x,x) + (x \ast y, x \ast y)$ and
$\partial_3((x,y,y))= - (x \ast y, y) + (x \ast y, y)$ for all $x, y \in X$.
If $x \sim y$, then $[x,x]= [y,y]$ in $H_2^{\rm D}(X)$.
Therefore $H_2^{\rm D}(X)$ is
generated by $\{ [x_\omega, x_\omega] | \omega \in {\rm Orb}(X) \}$,
where $x_\omega$ is a representative of an orbit $\omega$.
Since $H_2^{\rm D}(T_m)$ is
the free abelian group generated by $\{ [\omega, \omega] | \omega \in {\rm
Orb}(X)=T_m
\}$, we see that
$\pi_\ast: H_2^{\rm D}(X) \to H_2^{\rm D}(T_m)$ is an isomorphism.
Thus $H_2^{\rm D}(X)$ is a free abelian group of rank $m =|{\rm Orb}(X)|$.
In general,
from the diagram (\ref{eqn:BasicTrivializingDiagram}), we see that ${\rm
Ker}[i_\ast:
H_n^{\rm D}(X) \to H_n^{\rm R}(X)]$ is contained in
${\rm Ker}[\pi_\ast: H_n^{\rm D}(X) \to H_n^{\rm D}(T_m)]$.  Therefore we
have that
$i_\ast: H_2^{\rm D}(X) \to H_2^{\rm R}(X)$ is injective and hence
$\partial_\ast: H_3^{\rm
Q}(X) \to H_2^{\rm D}(X)$ is the $0$-map. Since
$\partial_\ast: H_2^{\rm Q}(X) \to H_1^{\rm D}(X)$ is the $0$-map, we have
the short exact sequence.
$\Box$

\begin{sect}{\bf Example.\/} {\rm
Let $X=R_k$ be the dihedral quandle of $k$ elements. Suppose that $k$ is an
odd integer.
Then $|{\rm Orb}(R_k)|=1$. Thus
\begin{eqnarray}
\begin{array}{ccccccccc}
0 & {\to}
& H_2^{\rm D}(X) & \stackrel{i_\ast}{\to} & H_2^{\rm R}(X) &
\stackrel{j_\ast}{\to} &
H_2^{\rm Q}(X) & {\to} & 0 \\
{} & {} & \downarrow \cong & {} & \downarrow & {} & \downarrow & {} & {} \\
0 & \to & H_2^{\rm D}(T_1) & \stackrel{i_\ast}{\to} &
H_2^{\rm R}(T_1) & \stackrel{j_\ast}{\to} & H_2^{\rm Q}(T_1) & \to & 0 \\
{} & {} & \| & {}
& \| & {} & \| & {} & {} \\
{} & {} & {\bf Z} & = & {\bf Z} &
{} & 0 & {} & {}
\end{array}
\end{eqnarray}
Greene proved that $H_2^{\rm R}(X)$ is generated by $[(0,0)]$ by a
geometric argument,
and the order is infinite by using the $\vec{\omega}$-orbit where $\vec{\omega}$ is the
generator of $H_2^{\rm
R}(T_1)$. Hence we have that
$H_2^{\rm Q}(X) = 0$. Conversely if we know that $H_2^{\rm Q}(X) = 0$, then
we have $H_2^{\rm R}(X) = {\bf Z}$. }\end{sect}

\begin{sect}{\bf Conjecture.\/}
In the basic homology long exact sequence for any finite quandle $X$, the
boundary operators
$\partial_\ast: H_n^{\rm Q}(X) \to
H_{n-1}^{\rm D}(X)$ are $0$-maps. Thus the sequence is decomposed into
short exacts
\begin{eqnarray}
0 \to H_n^{\rm D}(X)  \stackrel{i_\ast}{\to} H_n^{\rm R}(X)
\stackrel{j_\ast}{\to}
H_n^{\rm Q}(X) \to 0. \end{eqnarray}
\end{sect}

We define an index $S(X)$ of $X$
by the minimum integer $n$ such that
$\partial_\ast: H_n^{\rm Q}(X) \to
H_{n-1}^{\rm D}(X)$ is not the $0$-map (if there exist no such integers $n$,
then $S(X)=\infty$). The conjecture is that $S(X) = \infty$ for any finite
quandle $X$.

By a computer calculation, we have that $S(R_3)> 6$, $S(R_4)> 5$, $S(R_5)>4$,
$S(QS(5))>4$, etc. where $QS(5)$ is the quandle of non-identity permutations 
on three letters.

\section{Lower Bounds for Betti Numbers}

\begin{sect}\label{thm:BettiEstimB}
{\bf Theorem.\/}
Let $\pi: X \to {\rm Orb}(X) = T_m$ be the projection from a 
quandle $X$
to its orbit quandle. 
If  $X$ is finite  or 
if there is a homomorphism $s : T_m \to X$ with $\pi \circ s = {\rm id}$,
Then
\begin{eqnarray}
\beta_n^{\rm D}(X) \geq a_n, \quad
\beta_n^{\rm R}(X) \geq m^n, \quad
\beta_n^{\rm Q}(X) \geq b_n,
\end{eqnarray}
where $a_n$ and $b_n$ are  as before.
\end{sect}

Before proving this theorem, we give some remarks here.

(1) The inequalities of the theorem are best possible; namely, for any $n$,
there is a quandle $X$ such that the equalities hold. Actually, 
the
trivial
quandle $T_m$ is
such an example.

(2) In case $\pi_\ast$ is not surjective, the cokernel has a meaning.
Our proof of the theorem gives 
information on the cokernel that will
be treated later.

(3) Consider $X= {\bf R} \times T_2$ as a quandle with  
$$(a,i)*(b,j) =  \left\{ \begin{array}{lr} (2b-a, i) & {\rm if } \; i=j \\
                                          (a, i)     & {\rm if} \; i \ne j \end{array} \right.$$
Then this is a quandle with ${\rm Orb}(X) = T_2$ and with $s(j)=(0,j)$. In this case
$\pi \circ s = {\rm id}$. So the theorem applies to this infinite quandle.

Let $X$ be a finite quandle, 
and let $X \to {\rm Orb}(X)=T_m$ be the projection. 
For an $n$-tuple $\vec{\omega}= (\omega_1, \dots, \omega_n)$ 
of elements of ${\rm Orb}(X)$,
define an element $T^{\rm R}(\vec{\omega}) \in C_n^{\rm R}(X)$ by
\begin{eqnarray}
T^R(\vec{\omega}) =
\sum_{x_j \in \omega_j (j=1,\dots, n)} (x_1, \dots, x_n),
\end{eqnarray}
where $x_j$ runs over $\omega_j$
for each $j = 1, \ldots, n$.

For an $n$-tuple 
$\vec{\omega}= (\omega_1, \dots, \omega_n)$ of elements of ${\rm Orb}(X)$
such that $\omega_i =\omega_{i+1}$ for some
 $i \in \{1, \ldots , n-1\}$,
 pick an index $i_0 $
 such that $\omega_{i_0} =\omega_{i_0+1}$, and 
define an element 
$T^{\rm D}(\vec{\omega}; i_0) \in C_n^{\rm D}(X)$ by
\begin{eqnarray}
T^{\rm D}(\vec{\omega}; i_0)
 = \sum_{x_j \in \omega_j (j=1,\dots, n), x_{i_0}=x_{i_0+1}}
(x_1, \dots, x_n) \end{eqnarray}
where
$x_j$ runs over $\omega_j$ for each $j$ ($j=1, \dots, n$) under the
condition $x_{i_0} = x_{i_0+1}$.

\begin{sect} \label{lem:TransferA}
{\bf Lemma.\/}
\begin{itemize}
\item[{\rm (1)}]
$T^{\rm R}(\vec{\omega}) \in Z_n^{\rm R}(X)$.
\item[{\rm (2)}]
$T^{\rm D}(\vec{\omega}, i_0) \in Z_n^{\rm D}(X)$.
\end{itemize}
\end{sect}
{\it Proof.\/}
(1) If $n=1$, it is obvious. So we assume $n \geq 2$. \begin{eqnarray}
\lefteqn{
\partial_n(T^{\rm R}(\vec{\omega})) } \nonumber \\
&& =
\sum_{x_j \in \omega_j (j=1,\dots, n)} \partial_n(x_1, \dots, x_n) \nonumber \\
&& =
\sum_{x_j \in \omega_j (j=1,\dots, n)}
  \left[ \sum_{i=2}^{n} (-1)^{i}
  \left[ (x_1, x_2, \dots, x_{i-1}, x_{i+1},\dots, x_n) 
\right. \right. 
\nonumber \\
  && \phantom{aaaaaaaaaaaaaa}
  - 
\left. \left. 
(x_1 \ast x_i, x_2 \ast x_i, \dots,
  x_{i-1}\ast x_i, x_{i+1}, \dots, x_n) 
                             \right] \rule{0in}{20pt} \right]
\nonumber \\
&& =
\sum_{i=2}^{n} (-1)^{i}
\left[
\sum_{x_j \in \omega_j (j=1,\dots, n)}
\left[
(x_1, x_2, \dots, x_{i-1}, x_{i+1},\dots, x_n) \right.  \right. 
\nonumber \\
  && \phantom{aaaaaaaaaaaaaa}
  - 
\left. \left.
(x_1 \ast x_i, x_2 \ast x_i, \dots,
  x_{i-1}\ast x_i, x_{i+1}, \dots, x_n) 
\right] \rule{0in}{20pt}\right]
\end{eqnarray} Since
$S(x_i)|_{\omega_j}: \omega_j \to \omega_j$ is a bijection, there is a bijection
between the sets
\begin{eqnarray}
\{ (x_1, x_2, \dots, x_{i-1}, x_{i+1},\dots, x_n) |
x_j \in \omega_j (j=1,\dots, n) \} \nonumber
\end{eqnarray}
and
\begin{eqnarray}
\{ (x_1 \ast x_i, x_2 \ast x_i, \dots,
x_{i-1}\ast x_i, x_{i+1}, \dots, x_n) |
x_j \in \omega_j (j=1,\dots, n) \}. \nonumber
\end{eqnarray}
Thus the sum is zero.

(2) is proved by the same calculation.
$\Box$

\begin{sect} {\bf Proof of Theorem \ref{thm:BettiEstimB}.}
{\rm 
In the first case,
the induced homomorphism $s_\ast : H_n^{\rm W}(T_m) \to H_n^{\rm W}(X)$ is the
right inverse of $\pi_\ast : H_n^{\rm W}(X) \to H_n^{\rm W}(T_m)$. By
(\ref{eqn:BettiTrivial}), we have the inequalities.

In the second case,
 $H_n^{\rm R}(T_m) = C_n^{\rm R}(T_m)$ is a free abelian group generated by
 the $n$-tuples $\vec{\omega} =(\omega_1, \dots, \omega_n)$. Divide the generator set, ${\cal
G}_{\rm R}$, of
$C_n^{\rm R}(T_m)$ into two subsets ${\cal G}_{\rm D}$ and ${\cal G}_{\rm
Q}$ as follows:
${\cal G}_{\rm D}$
consists of $n$-tuples $\vec{\omega}=(\omega_1, \dots, \omega_n)$ such that
$\omega_i = \omega_{i+1}$ for some $i$, and ${\cal G}_{\rm Q}$ is the
complement.
For each generator 
$\vec{\omega} \in  
{\cal G}_{\rm D}$,
 fix an element
$T^{\rm D}(\vec{\omega}, i_0) \in Z_n^{\rm R}(X);$
for each
generator 
$\vec{\omega} \in {\cal G}_{\rm Q}$, 
consider the 
element
$T^{\rm R}(\vec{\omega}) \in 
Z_n^{\rm
R}(X)$. 
Obviously, 
$\pi_\ast(T^{\rm D}(\vec{\omega}, i_0))
= (\prod_{j=1}^{n}|\omega_j|)/|\omega_i| \vec{\omega}$, and
$\pi_\ast(T^{\rm R}(\vec{\omega})) =
(\prod_{j=1}^{n}|\omega_j|) \vec{\omega}$.
Define a homomorphism $T: H_n^{\rm R}(T_m) \to H_m^{\rm R}(X)$
by 
$$T(\vec{\omega})= \left\{ \begin{array}{lr} 
T^{\rm D}(\vec{\omega}, i_0) & {\rm if} \; \vec{\omega} \in {\cal G}_{\rm D}, \\
T^{\rm R}(\vec{\omega}) & {\rm if} \; \vec{\omega} \in {\cal G}_{\rm R}. \end{array} \right.$$
Then $\pi_\ast \circ T:
H_n^{\rm R}(T_m) \to H_n^{\rm R}(T_m)$ maps each generator in ${\cal
G}_{\rm D}$ (resp.
${\cal G}_{\rm Q}$) to itself multiplied by
$(\prod_{j=1}^{n}|\omega_j|)/|\omega_i|$ (resp.
$\prod_{j=1}^{n}|\omega_j|$). Thus the image of $T$ is a free abelian group
in $H_n^{\rm
R}(X)$ of rank $m^n$. Thus we have $\beta_n^{\rm R}(X) \geq m^n$.

$H_n^{\rm D}(T_m)$ is a subgroup of
$H_n^{\rm R}(T_m)$ generated by ${\cal G}_{\rm D}$. By
Lemma~\ref{lem:TransferA},
the image of the restriction of $T$ to $H_n^{\rm D}(T_m)$ is contained
in $H_n^{\rm D}(X)$, which
is a free abelian group of rank
$a_n$. Thus we have $\beta_n^{\rm D}(X) \geq a_n$.

The image of the restriction of $T$ to the subgroup of $H_n^{\rm R}(T_m)$
generated
by ${\cal G}_{\rm Q}$ is a free abelian group in $H_n^{\rm R}(X)$ of rank
$b_n$. 
Since
the subgroup of $H_n^{\rm R}(T_m)$ generated by ${\cal G}_{\rm Q}$ is
mapped identically to
$H_n^{\rm Q}(T_m)$, there is a free abelian group in $H_n^{\rm Q}(X)$ of
rank $b_n$. 
Thus we have
$\beta_n^{\rm Q}(X) \geq b_n$. 
$\Box$ }\end{sect}

\begin{sect}{\bf Corollary.\/}
Let $R_k$ be the dihedral quandle of $k$ elements. \begin{itemize}
\item[{\rm (1)}]
If $k$ is even, then
\begin{eqnarray}
\beta_n^{\rm D}(R_k) \geq 2^n-2, \quad
\beta_n^{\rm R}(R_k) \geq 2^n, \quad
\beta_n^{\rm Q}(R_k) \geq 2.
\end{eqnarray}
In particular,
$H_n^{\rm Q}(R_k;G) \neq 0$ and $H^n_{\rm Q}(R_k;G) \neq 0$ for any
coefficient group $G$.
\item[{\rm (2)}]
If $k$ is odd, then
\begin{eqnarray}
\beta_n^{\rm D}(R_k) \geq 1, \quad
\beta_n^{\rm R}(R_k) \geq 1. \quad
\end{eqnarray}
\end{itemize}
\end{sect}
{\it Proof.\/}
If $k$ is even, then $|{\rm Orb}(R_k)| = 2$. If $k$ is odd, then $|{\rm
Orb}(R_k)| = 1$.
By Theorem~\ref{thm:BettiEstimB}
and the universal coefficient theorem, we have the result. $\Box$

By a computer calculation, we have
\begin{eqnarray}
\beta_2^{\rm D}(R_4) = 2, \quad
\beta_2^{\rm R}(R_4) = 4, \quad
\beta_2^{\rm Q}(R_4) = 2.
\end{eqnarray}
\begin{eqnarray}
\beta_3^{\rm D}(R_4) = 6, \quad
\beta_3^{\rm R}(R_4) = 8, \quad
\beta_3^{\rm Q}(R_4) = 2.
\end{eqnarray}
Thus the lower bounds in the corollary (or Theorem~\ref{thm:BettiEstimB})
are the best possible.

\section{The Cokernel of $\pi_\ast$} 

Suppose the quandle $X$ is finite. 
Let $X \to {\rm Orb}(X)=T_m$ be the projection.
For an $n$-tuple 
$\vec{\omega}= (\omega_1, \dots, \omega_n)$ where $\omega_j \in {\rm Orb}(X),$ fix a
representative $x_n \in \omega_n$. 
Define 
\begin{eqnarray} 
T^{\rm R}(\vec{\omega}, x_n) =
\sum_{x_j \in \omega_j (j=1,\dots, n-1)} (x_1, \dots, x_n). \end{eqnarray}
The sum runs over
$x_j \in \omega_j$ for each $j=1, \dots, n-1$. 
Then
$T(\vec{\omega}, x_n)$
is an element of $C_n^{\rm R}(X)$.

Similarly, when
$\vec{\omega} =(\omega_1, \dots, \omega_n)$
and $\omega_j \in 
{\rm Orb}(X)$  
is 
such that
$\omega_i =\omega_{i+1}$ for some $i \in \{1,\dots, n-1\},$
we fix
representative $x_n \in
\omega_n$. 
Define 
\begin{eqnarray} 
T^{\rm D}(\vec{\omega}, i_0, x_n) =
\sum_{x_j \in \omega_j (j=1,\dots, n-1),\ x_{i_0}=x_{i_0+1}} (x_1, \dots, x_n) \end{eqnarray}
where  the sum runs over $x_j \in \omega_j$ such that 
$x_{i_0}=x_{i_0+1}$ and $x_n$ is fixed.
Then
$T^{\rm D}(\vec{\omega}, i_0, x_n)$
is an element of $C_n^{\rm D}(X)$.

By the same argument as in the proof of Lemma~\ref{lem:TransferA}, we see that
\begin{itemize}
\item[{\rm (1)}]
$T^{\rm R}(\vec{\omega}, x_n) \in Z_n^{\rm R}(X)$,
\item[{\rm (2)}]
$T^{\rm D}(\vec{\omega}, i_0, x_n) \in Z_n^{\rm D}(X)$.
\end{itemize}
In the proof of Theorem~\ref{thm:BettiEstimB}, we may consider
a homomorphism  $T': H_n^{\rm R}(T_m) \to H_n^{\rm R}(X)$,
instead of $T$, such that
$$T'(\vec{\omega})= \left\{ \begin{array}{lr} 
T^{\rm D}(\vec{\omega}, i_0, x_n) & {\rm if} \; \vec{\omega} \in {\cal G}_{\rm D} \\
T^{\rm R}(\vec{\omega}, x_n) & {\rm if} \;  \vec{\omega} \in {\cal G}_{\rm R} \end{array} \right.$$
Then $\pi_\ast \circ T':
H_n^{\rm R}(T_m) \to H_n^{\rm R}(T_m)$ maps each generator
in ${\cal G}_{\rm D}$ (resp. ${\cal G}_{\rm Q}$) to
itself multiplied 
by 
$(\prod_{j=1}^{n-1}|\omega_j|)/|\omega_i|$ (resp.
$\prod_{j=1}^{n-1}|\omega_j|$).
   The cokernel of $\pi_\ast: H_n^{\rm W}(X) \to H_n^{\rm W}(T_m)$ is
generated by
${\cal G}_{\rm D}$, ${\cal G}_{\rm R} = {\cal G}_{\rm D}\cup {\cal G}_{\rm Q}$,
or ${\cal G}_{\rm Q}$, according to ${\rm W}$ is ${\rm
D}$, ${\rm R}$ or ${\rm Q}$.
If 
$\vec{\omega}$
is in ${\cal G}_{\rm D}$, its
order in
${\rm Coker}[\pi_\ast: H_n^{\rm W}(X) \to H_n^{\rm W}(T_m)]$
(${\rm W} = {\rm D}, {\rm R}$) is finite and divides
$(\prod_{j=1}^{n-1}|\omega_j|)/|\omega_i|$.
If 
$\vec{\omega}$
is in ${\cal G}_{\rm Q}$, its order in
${\rm Coker}[\pi_\ast: H_n^{\rm W}(X) \to H_n^{\rm W}(T_m)]$
(${\rm W} = {\rm R}, {\rm Q}$) is finite and divides
$\prod_{j=1}^{n-1}|\omega_i|$.
Here we assume that 
a trivial element has order $1$. Therefore we have the
following.

\begin{sect}\label{thm:BettiEstimC}
{\bf Proposition.\/}
Let $\pi: X \to {\rm Orb}(X) = T_m$ be the projection from a finite quandle $X$
to its orbit quandle. The cokernel of $\pi_\ast: H_n^{\rm W}(X) \to
H_n^{\rm W}(T_m)$ is
finite. The order of each generator
$[\vec{\omega}]$ in
${\rm Coker}[\pi_\ast: H_n^{\rm W}(X) \to H_n^{\rm W}(T_m)]$ divides
$(\prod_{j=1}^{n-1}|\omega_j|)/|\omega_i|$
if $\vec{\omega} \in {\cal G}_{\rm D}$ or
$\prod_{j=1}^{n-1}|\omega_i|$ if
$\vec{\omega} \in {\cal G}_{\rm Q}$. \end{sect}

We abbreviate ${\rm Coker}[\pi_\ast: H_n^{\rm W}(X) \to H_n^{\rm W}(T_m)]$ to
${\rm Coker}_n^{\rm W}(X)$.
For $R_4$, $|{\rm Orb}(0)| = |{\rm Orb}(1)|= 2$. If $n=2$, then
${\cal G}_{\rm D}= \{(\omega_0, \omega_0), (\omega_1, \omega_1)\}$ and
${\cal G}_{\rm Q}= \{(\omega_0, \omega_1), (\omega_1, \omega_0)\}$, where
$\omega_i= {\rm
Orb}(i)$. By Proposition~\ref{thm:BettiEstimC}, we have that ${\rm
order}(\omega_0,
\omega_0) = {\rm order}(\omega_1, \omega_1) = 1$
and that
${\rm order}(\omega_1, \omega_0)$ and
${\rm order}(\omega_0, \omega_1)$ are $1$ or $2$. Thus
\begin{eqnarray}
{\rm Coker}_2^{\rm D}(R_4) = 0, \quad
{\rm Coker}_2^{\rm R}(R_4) = ({\bf Z}_2)^k, \quad
{\rm Coker}_2^{\rm Q}(R_4) = ({\bf Z}_2)^k, \end{eqnarray}
for some $k \in \{0,1,2\}$.
By a computer calculation, we have
\begin{eqnarray}
{\rm Coker}_2^{\rm D}(R_4) = 0, \quad
{\rm Coker}_2^{\rm R}(R_4) = ({\bf Z}_2)^2, \quad
{\rm Coker}_2^{\rm Q}(R_4) = ({\bf Z}_2)^2. \end{eqnarray}

\section{Another Relation between the Homology Groups}

We will
give an alternative relationship between the degeneration homology groups and
the rack homology groups.

For a quandle $X$, let $C_n^{\rm DD}(X)$ be the free abelian group generated
by $n$-tuples $(x_1, \dots, x_n)$ such that $x_1=x_2$, or
$C_n^{\rm DD}(X)=0$ if $n <2$.  Then $C_\ast^{\rm DD}(X) =\{ C_n^{\rm DD}(X),
\partial_n\}$  is a sub-complex of $C_\ast^{\rm D}(X)$ and of $C_\ast^{\rm
R}(X)$.
Putting $C_n^{\rm D/DD}(X) = C_n^{\rm D}(X) / C_n^{\rm DD}(X)$, we have
a chain complex $C_\ast^{\rm D/DD}(X) =\{ C_n^{\rm D/DD}(X), \partial_n\}$
and
a long exact sequence
\begin{eqnarray}\label{eqn:HomologyLongExactD}
\cdots \stackrel{\partial_\ast}{\to} H_{n}^{\rm DD}(X;G)
\stackrel{i_\ast}{\to} H_n^{\rm D}(X;G)
\stackrel{j_\ast}{\to} H_n^{\rm D/DD}(X;G)
\stackrel{\partial_\ast}{\to} H_{n-1}^{\rm DD}(X;G) \to \cdots
\end{eqnarray}

\begin{sect}\label{thm:HomologyLongExactRD}
{\bf Proposition.\/}
For a quandle $X$, there is
a long exact sequence
\begin{eqnarray}\label{eqn:HomologyLongExactRD}
\cdots {\to} H_{n-1}^{\rm R}(X;G)
{\to} H_n^{\rm D}(X;G)
\stackrel{j_\ast}{\to} H_n^{\rm D/DD}(X;G)
{\to} H_{n-2}^{\rm R}(X;G) \to \cdots
\end{eqnarray}
This is natural with respect to homomorphisms induced from quandle
homomorphisms.
\end{sect}
\begin{sloppypar}
{\it Proof.\/}
Let $u_n: C_n^{\rm DD}(X) \to C_{n-1}^{\rm R}(X)$ be an isomorphism
with $u_n(x_1, \dots, x_n)$ $=$ $(x_1, x_3, \dots, x_n)$.
It is easily checked that $u_{n-1} \circ \partial_n = - \partial_{n-1}
\circ u_n$,
namely,
$u= \{ u_n \}: C_\ast^{\rm DD}(X) \to C_\ast^{\rm R}(X)$
is a chain map of degree $-1$.
It induces an isomorphism ${u_n}_\ast : H_n^{\rm DD}(X) \to H_{n-1}^{\rm R}(X)$.
Combine this isomorphism with (\ref{eqn:HomologyLongExactD}).
$\Box$\end{sloppypar}

For a quandle $X$, let $UH_n^{\rm DD}(X)$ and $UH_n^{\rm D}(X)$ be
the subgroups of $H_n^{\rm DD}(X)$ and $H_n^{\rm D}(X)$ generated
by  $\{ [(x, \dots, x)] \/ | \/ x \in X \}$.

\begin{sect}\label{lem:kerDD}
{\bf Lemma.\/}
${\rm Ker}[ i_\ast : H_n^{\rm DD}(X) \to H_n^{\rm D}(X) ]
\cap UH_n^{\rm DD}(X) = 0.$
\end{sect}
{\it Proof.\/}
Let $\pi: X \to {\rm Orb}(X) = T_m$ be the projection
to its orbit quandle.
{}From the naturality of the exact sequence~(\ref{eqn:HomologyLongExactD}),
we have a commutative diagram
\begin{eqnarray} 
\begin{array}{ccccccccc}
\cdots & \stackrel{\partial_\ast}{\to}
& H_n^{\rm DD}(X) & \stackrel{i_\ast}{\to}
& H_n^{\rm D}(X) & \stackrel{j_\ast}{\to}
& H_n^{\rm D/DD}(X) &
\stackrel{\partial_\ast}{\to} & \cdots
\\ {} & {} & \downarrow & {} & \downarrow & {} & \downarrow & {} & {} \\
0 & \to
& H_n^{\rm DD}(T_m) & \stackrel{i_\ast}{\to}
& H_n^{\rm D}(T_m) & \stackrel{j_\ast}{\to}
& H_n^{\rm D/DD}(T_m) & \to & 0, \end{array}
\end{eqnarray}
where the vertical maps are the induced homomorphisms $\pi_\ast$.
By a similar argument to the proof of Lemma~\ref{lem:SecondHomology},
we see that $UH_n^{\rm DD}(X)$ and $UH_n^{\rm D}(X)$ are
generated by $\{ [x_\omega, \dots, x_\omega] | \omega \in {\rm Orb}(X) \}$,
where $x_\omega$ is a representative of an orbit $\omega$ and that
$\pi_\ast: UH_n^{\rm DD}(X) \to UH_n^{\rm DD}(T_m)$
and
$\pi_\ast: UH_n^{\rm D}(X) \to UH_n^{\rm D}(T_m)$ are isomorphisms.
Note that $UH_n^{\rm DD}(T_m)$ and $UH_n^{\rm D}(T_m)$ are
free abelian group generated by
$\{ [\omega, \dots, \omega] | \omega \in {\rm Orb}(X) \}$.
Thus we have the result. $\Box$

\begin{sect}\label{lem:DD3andDD4}
{\bf Lemma.\/}
The boundary operators
$\partial_\ast: H_4^{\rm D/DD}(X) \to H_3^{\rm DD}(X)$ and
$\partial_\ast: H_3^{\rm D/DD}(X) \to H_2^{\rm DD}(X)$ are
$0$-maps.
\end{sect}
{\it Proof.\/}
Let $s_n: C_n^{\rm D/DD}(X) \to C_n^{\rm D}(X)$ be a homomorphism
defined by $(x_1, \dots, x_n) + C_n^{\rm DD}(X) \mapsto (x_1, \dots, x_n)$
where $(x_1, \dots, x_n)$ are $n$-tuples such that $x_1 \neq x_2$
and there exists some $i$ with $x_i = x_{i+1}$.
There is a unique homomorphism
$\phi_n : C_n^{\rm D/DD}(X) \to C_{n-1}^{\rm DD}(X)$ such that
$i \circ \phi_n = \partial_n \circ s_n - s_{n-1} \circ \partial_n$.
Then $\phi =\{\phi_n\}$ is a chain map of degree $-1$, i.e.,
$\partial_{n-1} \circ \phi_n = - \phi_{n-1} \circ \partial_n$, and
the induced homomorphism
$(\phi_n)_\ast: H_n^{\rm D/DD}(X) \to H_{n-1}^{\rm DD}(X)$
is the same as the boundary operator
$\partial_\ast: H_n^{\rm D/DD}(X) \to H_{n-1}^{\rm DD}(X)$.
For simplifying notation, we denote an element
$(x_1, \dots, x_n) + C_n^{\rm DD}(X)$ of $C_n^{\rm D/DD}(X)$ by
$(x_1, \dots, x_n)$.

$C_4^{\rm D/DD}(X)$ is generated by $(x_1,x_2,x_2,x_4)$
and $(x_1,x_2,x_3,x_3)$ for $x_1,\dots, x_4 \in X$ with
$x_1 \neq x_2$.  Since
$\phi_4((x_1,x_2,x_2,x_4))= 0$
and
$\phi_4((x_1,x_2,x_3,x_3))= 0$ or $\pm (x_3,x_3,x_3)$,
the image ${\rm Im}[\partial_\ast: H_4^{\rm D/DD}(X) \to H_3^{\rm D}(X)] $
is in $UH_3^{\rm DD}(X)$.  By Lemma~\ref{lem:kerDD} and the
exactness 
of (\ref{eqn:HomologyLongExactD}), we have
$${\rm Im}[\partial_\ast: H_4^{\rm D/DD}(X) \to H_3^{\rm D}(X)] = 0.$$

Since $H_2^{\rm DD}(X) = UH_2^{\rm DD}(X)$, by Lemma~\ref{lem:kerDD}
and the exactness of (\ref{eqn:HomologyLongExactD}),
$${\rm Im}[\partial_\ast: H_3^{\rm D/DD}(X) \to H_2^{\rm D}(X)] = 0.
\Box $$ 

\begin{sect}\label{lem:DD3}
{\bf Lemma.\/}
$\pi_\ast : H_3^{\rm D/DD}(X) \to H_3^{\rm D/DD}(T_m)$ is an isomorphism
where $T_m = {\mbox{ \rm Orb}}(X)$.
\end{sect}
\begin{sloppypar}
{\it Proof.\/}
$C_3^{\rm D/DD}(X)$ is generated by $(x_1,x_2,x_2)$
for $x_1, x_2 \in X$ with $x_1 \neq x_2$.  Since
$\partial_3( (x_1,x_2,x_2) ) =0, $ 
we have that
$Z_3^{\rm D/DD}(X) = C_3^{\rm D/DD}(X)$ and $H_3^{\rm D/DD}(X)$ is
generated by $[(x_1,x_2,x_2)]$ ($x_1 \neq x_2$).
Since
$\partial_4(x_1, x_2, x_2, x_4) = (x_1, x_2, x_2) - (x_1 \ast x_4,
x_2 \ast x_4, x_2 \ast x_4)$
and
$\partial_4(x_1, x_2, x_3, x_3) = (x_1, x_3, x_3) - (x_1 \ast x_2,
x_3, x_3)$,
we see that if $x_1 \sim x'_1$ and $x_2 \sim x'_2$, then
$[(x_1, x_2, x_2)] = [(x'_1, x'_2, x'_2)]$ in $H_3^{\rm D/DD}(X)$.
If $x_1 \sim x_2$, then
$[(x_1, x_2, x_2)] = [(x_1, x_1, x_1)]= 0$ in $H_3^{\rm D/DD}(X)$.
Thus $H_3^{\rm D/DD}(X)$ is
generated by $[(x_{\omega_1},x_{\omega_2},x_{\omega_2})]$
($\omega_1, \omega_2 \in {\rm Orb}(X)$ with $\omega_1 \neq \omega_2$),
where $x_\omega$
is a representative of $\omega \in {\rm Orb}(X)$.
Since $H_3^{\rm D/DD}(T_m)$ is a free abelian group generated by
$[(\omega_1,\omega_2,\omega_2)]$
($\omega_1, \omega_2 \in {\rm Orb}(X)$ with $\omega_1 \neq \omega_2$),
we have the result. $\Box$ \end{sloppypar}

\begin{sect}
{\bf Proposition.\/}
For a quandle $X$, there exists a short exact sequence
\begin{eqnarray}
0 \to H_2^{\rm R}(X)  {\to} H_3^{\rm D}(X)  {\to}
{\bf Z}^{m^2 -m} \to 0,  \end{eqnarray}
where $m = |{\rm Orb}(X)|$.
\end{sect}
{\it Proof.\/}
By Lemma~\ref{lem:DD3andDD4}, we have a short exact sequence
\begin{eqnarray}
0 \to H_2^{\rm R}(X)  {\to} H_3^{\rm D}(X)  {\to}
H_3^{D/DD}(X) \to 0   \end{eqnarray}
from (\ref{eqn:HomologyLongExactRD}).
By Lemma~\ref{lem:DD3}, $H_3^{D/DD}(X)$ is isomorphic to ${\bf Z}^{m^2
-m}$.  $\Box$

\begin{sect}
{\bf Corollary.\/}
${\rm Torsion}H_2^{\rm R}(X) \cong {\rm Torsion}H_3^{\rm D}(X).$
\end{sect}

\section{ Quandle (Co)homology and Virtual Knots} 

\begin{figure}
\begin{center}
\mbox{
\epsfxsize=3in 
\epsfbox{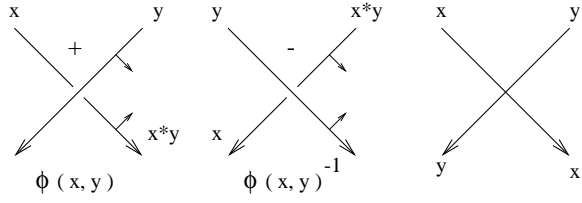}
}
\end{center}
\caption{Three types of crossings of virtual knots}
\label{threecrossings}
\end{figure}

A {\it  virtual knot (diagram)} \cite{Lou} is a
generically  immersed   oriented $1$-manifold
in the plane together with the following three types of  crossing information
at double points.
First, there are two types, positive and negative, crossings
with over-under information as in the classical knot theory.
The under-path is broken into two arcs. The left and the middle pictures
of Fig.~\ref{threecrossings}
depict positive and negative crossings, respectively.
(The labels and the expression $\phi$ will be used later.) 
The right of the figure depicts a crossing of the third type, 
called a {\it virtual crossing}, at which there is no over-under information.

\begin{figure}
\begin{center}
\mbox{
\epsfxsize=3in 
\epsfbox{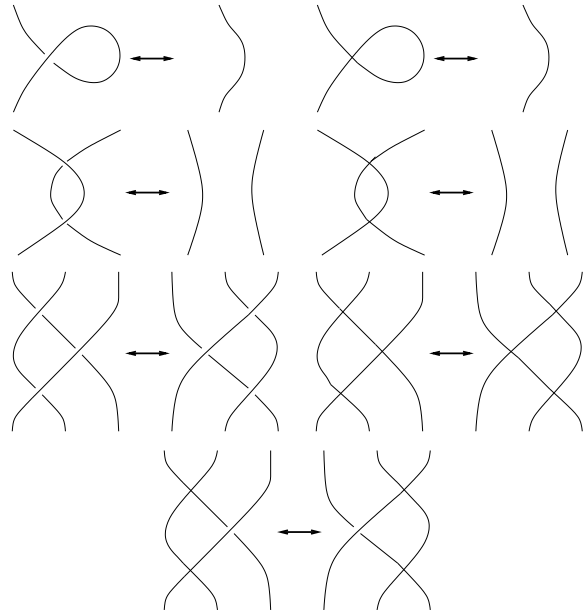}
}
\end{center}
\caption{Reidemeister moves for virtual knots by Kauffman}
\label{virIII}
\end{figure}

Two virtual knot diagrams  are {\it equivalent} if the diagrams 
are related by a sequence of 
Reidemeister moves depicted in Fig.~\ref{virIII}, 
and ambient isotopy of the plane. 
A {\it virtual knot} is an equivalence class of a virtual knot diagram.

\begin{figure}
\begin{center}
\mbox{
\epsfxsize=3in
\epsfbox{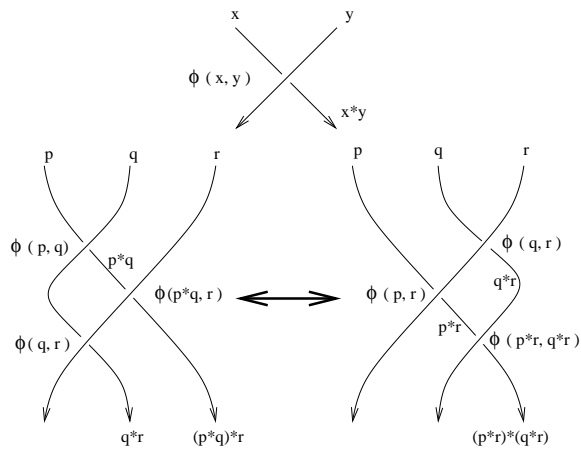}
}
\end{center}
\caption{Colors at a crossing }
\label{2cocy}
\end{figure}

At crossings of a virtual knot, the under-arc is broken.
The complement consists of immersed arcs.   
These transverse components of arcs are called over-arcs of a virtual knot.

A {\it color} (or {\it coloring}) 
on a virtual 
knot 
diagram is a
function ${\cal C} : R \rightarrow X$, where $X$ is a fixed 
quandle and $R$ is the set of
over-arcs 
  satisfying the  condition
depicted in the top
of Fig.~\ref{2cocy}. 
In  the top of 
Fig.~\ref{2cocy}, 
a crossing with over-arc, $r$, has color ${\cal C}(r)= y \in X$. 
The under-arcs are called $r_1$ and $r_2$ from top to bottom; they are colored 
${\cal C}(r_1)= x$ and ${\cal C}(r_2)=x*y$.
Note that locally the colors do not depend on the 
orientation of the under-arc.

Assume that a finite quandle $X$ is given. 
Pick a 
quandle 
$2$-cocycle 
$\phi \in  Z_{\rm Q}^2(X, G),$ 
and write the coefficient 
group, $G$, multiplicatively. 
Consider a non-virtual crossing 
in the diagram.
For each coloring of the diagram, evaluate 
the $2$-cocycle 
 on the quandle colors that 
appear near the crossing as described as follows: 
The first argument is the color on the under-arc away from which 
the normal to the over-arc points. 
The second argument is the color on the over-arc. See Fig.~\ref{2cocy}.

Let $\tau$ denote 
a non-virtual crossing, 
let $\epsilon (\tau)$ denote its sign, 
and let ${\cal  C}$ denote a coloring.
When the colors of 
the arcs 
 are as describe above,
the 
{\it (Boltzmann) weight
of a crossing} is  
$B(\tau, {\cal C}) = \phi(x,y)^{\epsilon (\tau)}$. 

The {\it partition function}, or a {\it state-sum}, 
is the expression 
$$
\sum_{{\cal C}}  \prod_{\tau}  B( \tau, {\cal C}).
$$
The product is taken over all crossings of the given diagram,
and the sum is taken over all possible colorings.
The values of the partition function 
are  taken to be in  the group ring ${\bf Z}[G]$ where $G$ is the coefficient 
group.
In fact, the value is in the group ``rig'' ${\bf N}[G]$.

By checking the equivalence relations we obtain

\begin{sect} {\bf Proposition.\/} 
 The state-sum is invariant 
under equivalence relations for 
virtual 
knots, 
thus defines invariants 
 $\Phi (K)$
(or $\Phi_{\phi}(K)$ to specify the 
$2$-cocycle $\phi $ used). 
\end{sect}

\begin{sect} {\bf Proposition.\/}  \label{coblemma}
If $\Phi_{\phi}$ and $\Phi_{\phi '} $ denote the state-sum invariants 
defined from cohomologous $2$-cocycles  $\phi$ and $\phi'$
then $\Phi_{\phi} =\Phi_{\phi '} $ (so that $\Phi_{\phi} (K)=\Phi_{\phi '}(K)$
 for any 
classical knot, or  virtual knot).  
In particular, 
the state-sum is equal to the number of
colorings of
a given 
knot diagram    
if the $2$-cocycle 
used for the Boltzmann weight is a coboundary.
\end{sect}

\begin{sect}{ \bf Remark.\ }{\rm 
 The definition of
colors and the above propositions naturally generalize those in
\cite{CJKLS}, stated for classical knots,
to virtual knots.  
The state-sum invariants are defined also for knotted surfaces in $4$-space 
in \cite{CJKLS} and studied in \cite{CJKS}.
For surfaces, $3$-cocycles are used as Boltzmann weights assigned to 
triple points on projections.
Virtual knots can also be defined in higher dimensions.
  A detailed study of these will be presented in a forthcoming paper.
}\end{sect}

We use the notion of linking numbers of virtual links \cite{GPV} in 
 the next section for construction of examples.
Let $L=K_1 \cup K_2$ be a virtual 
 link, where
$K_i$ ($i=1,2$) are distinct components.
Let $P$ and $N$ be the numbers of positive and negative, respectively, 
crossing of $L$ such that at the crossings $K_1$ goes over $K_2$.
Define the {\it virtual linking number} $\vlk(K_1, K_2)=P-N$.

\begin{figure}
\begin{center}
\mbox{
\epsfxsize=2.75in 
\epsfbox{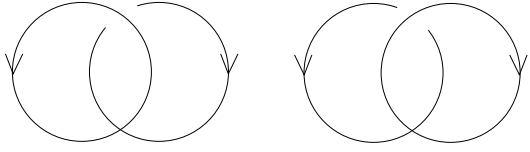}
}
\end{center}
\caption{Negative and positive
virtual Hopf links}
\label{hopf}
\end{figure}

\begin{sect} {\bf Lemma.\/}
 For any prescribed integers $n_{ij}$, $i,j=1, \cdots, k$,
there is a virtual link $L=K_1 \cup \cdots \cup K_k$ such that 
$\vlk(K_i,K_j)=n_{ij}$.
\end{sect}
{\it Proof.\/} 
Consider 
a virtual Hopf link $H_{\pm}=K_1 \cup K_2$,
 the Hopf link diagram
with one $\pm$-crossing respectively and one virtual crossing
(Fig.~\ref{hopf}).
If the first component goes over the second, $\vlk(K_1, K_2)=\pm1$
and $\vlk(K_2, K_1)=0$. 
The result follows by taking appropriate connected sum of copies of these.
$\Box$

\begin{sect} {\bf Proposition.\/} 
The cocycle invariants with trivial quandles $T_n$ depends only on 
the virtual linking numbers.
\end{sect}
{\it Proof.\/} The colors are constant on each component.
Any cocycle is written as a product of characteristic functions
$\chi_{(i,j)}$, so the state-sum is described by $\vlk$.
$\Box$

\section{Applications to Quandle (Co)homology}

Let $\pi: X \rightarrow Y$ be a surjective quandle homomorphism.
Since $Y$ is generally smaller, we try to use the information 
we already have for (co)homology groups of $Y$ to obtain new information
for those of $X$. Here, we apply this technique to a variety of quandles.
The coefficients of the (co)homology groups are ${\bf Z}$ unless 
otherwise specified.

\begin{figure}
\begin{center}
\mbox{
\epsfxsize=4in 
\epsfbox{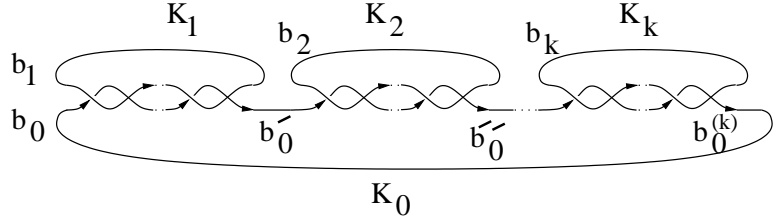}
}
\end{center}
\caption{A family of virtual links }
\label{virtori}
\end{figure}

\begin{sect}{\bf Proposition. \ } Let a  virtual knot or link  diagram $K$  be 
colored by a quandle $X$. 
Then $K$ represents a 2-cycle in $Z^{\rm W}_2(X)$ where ${\rm W}= {\rm R}$ or ${\rm Q}$.
\end{sect}
 {\it Proof.} Consider a (non-virtual) crossing  of $K$. Then the colors
 $(x,y)$ 
that are adjacent to the crossing represent a chain. 
As usual, $x$ is the color on the under-arc away from which
the normal to the over crossing points, and $y$ is the color on the over-arc.
We define the sign of such a chain to be the sign of the crossing.
The sum of these signed chains (taken over all crossings) is clearly a cycle.
 $\Box$

\begin{sect} {\bf Theorem.\/}   \label{nontrivLX2}
Let $X={\bf Z}[T,T^{-1}]/(h(T))$ be an Alexander quandle
 where $h$ is a polynomial with $h(1)=0$,
$T_{\infty}={\bf Z}$ be the trivial quandle, and 
$\pi: X \rightarrow T_{\infty}$ be the quandle homomorphism 
defined by $\pi(f(T))=f(1)$.
Then
the homomorphisms $\pi_* : H^{\rm Q}_2(X) \rightarrow H^{\rm Q}_2(T_{\infty})$
and $\pi^*: H_{\rm Q}^2(T_{\infty}) \rightarrow  H_{\rm Q}^2(X)$ are not $0$-maps.
In particular, $H^{\rm Q}_2(X) \neq 0 \neq H_{\rm Q}^2 (X)$.
\end{sect}
{\it Proof.\/}
For a given $X$ with $h(1)=0$, it will be proved in the lemma
that follows this proof that there is a virtual link $L$ 
(depicted in Fig.~\ref{virtori}, connected sums of virtual torus links)
such that (1) $L$ has a nontrivial color with $X$,  and (2) the color 
contributes a
nontrivial $t$-term to the state-sum with the cocycle 
$\pi^\sharp(\chi_{(a,b)})$ for some $a, b \in T_{\infty}$, 
where $\chi$ denotes the characteristic function:
$$\chi_x(y) = \left\{ \begin{array}{lr} 1 & {\mbox{\rm if}} \  \ x=y,
\\
                                               0   & {\mbox{\rm if}} \ \  x\ne
y. \end{array} \right.  $$  
Hence by Proposition~\ref{coblemma}, $\pi^\sharp(\chi_{(a,b)})$ is not a 
coboundary, and $\pi^*: H_{\rm Q}^2(T_{\infty}) \rightarrow H_{\rm Q}^2(X)$ 
is not the $0$-map.

The above color of $L$ by $X$ 
determines a $2$-cycle $\alpha$ in $Z^{\rm Q}_2(X)$ as in the preceding proposition. 
The $2$-cycle $\pi_\sharp(\alpha)$ is represented
by the same link $L$ with the colors taken mod $(T-1)$, 
i.e., the colors with substitution $T=1$. 
There are crossings in $L$ with different colors $(a, b)$,
 $a \neq b \in T_{\infty}$ after substitution $T=1$. Therefore 
 $\pi_\sharp(\alpha)$ 
is not a coboundary in $Z^{\rm Q}_2(T_\infty)$; so $\alpha$ is non-trivial in $H^{\rm Q}_2(X)$,
and $\pi_*: H^{\rm Q}_2(X) \rightarrow H^{\rm Q}_2( T_{\infty})$ is not the $0$-map.
$\Box$

\begin{sect} {\bf Lemma.\/}
The virtual link $L$ depicted in Fig.~\ref{virtori}
(where the numbers of crossings will be determined in the proof for 
any given $X$)
is colored nontrivially by $X$, and 
has the nontrivial state-sum term with this color.
\end{sect}
{\it Proof.\/}
 Let $h(T)=c_0+\sum_{i=1}^k c_i T^{m_i}$ be a
polynomial such that $ c_i \neq 0$ for $i=1,2,\ldots,n$ and
$\{ m_i \}_{i=1}^k$ is 
a 
strictly increasing sequence of positive
integers. Then any polynomial $h(T)$ with $h(1)=0$ 
can be written as such a polynomial if and only if $\sum_{i=0}^k c_i=0,$
which is equivalent to   $c_k= - \sum_{i=0}^{k-1} c_i.$ 
In Fig.~\ref{virtori}, 
the crossing repeats the sequence of
 a positive crossing followed by a virtual crossing, and for
$i=1,2,\ldots,k,$ let $\vlk(K_i,K_0)=n_i$, where $n_i$ will
be specified below, and all other virtual
linking numbers to be 0.
 Color each $K_i$ initially by $b_i.$ 
Note
that the color  $b_0$ changes 
to $b'_0,b''_0, \ldots
,b^{(k)}_0$ as the string $K_0$ links with the other components 
$K_1, \cdots, K_k$, as depicted in the figure.  
Now, let $n_k=m_1$ and for $i=1,2,\ldots, k-1,$
$n_{k-i}=m_{i+1}-m_i.$ We see that
$b^{(i)}_0=T^{n_i} b^{(i-1)}_0 + (1-T^{n_i}) b_i,$ so inductively,
$$b^{(k)}_0=b_k+\sum_{i=1}^{k}(b_{k-i}-b_{k-i+1}) 
T^{m_i}.$$
Take $b_{k-j}=\sum^{j}_{t=0}c_j$ for
$j=0,1,\ldots, k-1$ and $b_0=0$. 
{}From these definitions, we see
that $b_{k-i}-b_{k-i+1}=c_i$ 
for $i=1,2, \ldots, k-1,$ 
and  $-b_1=-\sum_{i=0}^{k-1}c_i=c_k.$
The right-hand-side of the expression for $b_0^{(k)}$ is 
$h(T)$.
Since $h(T)$ is equivalent to $0$ in 
$X$, we have a coloring ${\cal C}$ of $L$.
With the $2$-cocycle
$\pi^{\sharp}(\chi_{(0,b_1)})$,	 
the  state-sum term, $\prod_{\tau} B(\tau, {\cal C})$ is $t$
to the power
at least $\vlk(K_1,K_0),$ which is not an integer.
$\Box$

\begin{figure}
\begin{center}
\mbox{
\epsfxsize=2in 
\epsfbox{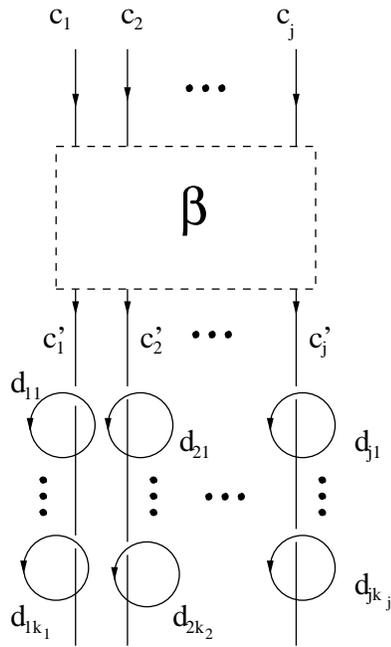}
}
\end{center}
\caption{Braids with virtual circles }
\label{vbraid}
\end{figure}

\begin{figure}
\begin{center}
\mbox{
\epsfxsize=4in 
\epsfbox{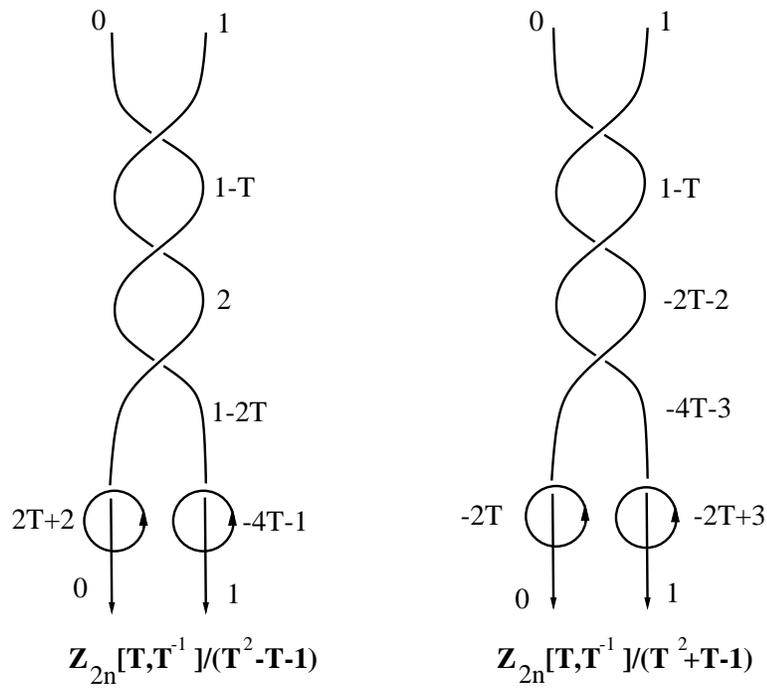}
}
\end{center}
\caption{Trefoils with virtual circles }
\label{trefoils}
\end{figure}

\begin{sect}{\bf Theorem.\/} Let $X$ and $Y$ be quandles.
  Suppose there exists $\pi : X \rightarrow Y$,
 a surjective homomorphism that is locally-homogeneous, and  there is a
link $L$ 
and a cocycle $\phi$  in $Z_{\rm Q}^2(Y, G)$
such that $\Phi_{\phi}(L)$ is non-trivial. Then $H_{\rm Q}^2(X, G) \neq
0.$
\end{sect}
{\it Proof.\/} 
Let $L$ be the link such that
$\Phi_{\phi}(L)$ is non-trivial (i.e., not an integer). 
 To prove $H_{\rm Q}^2(X, G) \neq 0,$ it is
sufficient to show that there is a (virtual) link $K$, a cocycle
$\psi \in Z_{\rm Q}^2(X, G),$ and a color ${\cal C}$ such that
 $\prod_{\tau} B (\tau,{\cal C}),$ the state-sum term for $K$ associated to
${\cal C},$ is not an integer.  To construct such a 
virtual link
$K$ and a 
color
${\cal C},$ first start with $L$ colored by $Y.$  Note that $K$
may be considered as the closed form of a $j$-strand braid,
$\beta,$ for some $j \in {\bf Z}.$   Since $\Phi_{\phi}(L)$ is
non-trivial, there exists a color ${\cal C'}$ of 
$L$ (regarded as a closed braid) 
such that the state-sum term associated to ${\cal C'}$ is
non-trivial. Observe that ${\cal C'}$ can be uniquely represented
as a choice of colors $b_1,b_2, \ldots ,b_j$ on the initial (top)
segments of $\beta.$  We now start constructing 
a virtual link
$K,$ and its color  ${\cal C}.$ Begin with the braid $\beta,$ and
color it initially (at the top) by $c_1,c_2, \ldots ,c_j,$ where
$c_i \in \pi^{-1}(b_i)$ for $i=1,2,\ldots,j,$ and extend the color to
all the segments of the braid $\beta$.
 Note that since
$\pi$ is a homomorphism, if a segment of $\beta$ is labeled $g$
when colored by $Y,$ the segment will be labeled $\pi(g)$ when
colored by $X.$  Thus, the terminal ends of $\beta,$  
are colored by  $c'_1,c'_2, \ldots ,c'_j,$ 
with
the property
that $\pi(c_i)=\pi(c'_i)$ for $i=1,2,\ldots,j.$  Since $\pi$ is
locally-homogeneous, there exists 
a word $w_i = d_{j_1}^{\epsilon_1} d_{j_2}^{\epsilon_2}
 \cdots d_{j_{k_i}}^{\epsilon_{k_i}}$
where each $d_{j_m} \in \pi^{-1}(b_i)$
such that $c'_i *w_i =c_i.$
For each strand $i$
of $\beta$ attach $k_i$ simple closed loops $K_{i_1}, K_{i_2},\ldots
K_{i_{k_i}}$ that cross over strand $i$ and returns via a virtual
crossing such that $\vlk(K_{i_m},K)=\epsilon_m,$ and $\vlk(K,K_{i_m})=0.$
See Fig.~\ref{vbraid}. Color
each $K_{i_m}$ by $d_{i_m}$. 
The closure of the
braid with the virtual loops is the virtual link $K$
we needed. 
  Take
$\psi=\pi^\sharp(\phi),$ and notice that the new crossings created by
the added virtual links  have trivial state-sum contributions.
Hence, the state-sum term 
of $K$ for ${\cal C}$ with respect to
$\psi$ is equal to
 the state-sum term of $K$ for ${\cal C'}$
 with respect to $\phi,$ and so is non-trivial.
$\Box$

Note that for a given link $L$ and a color ${\cal C}$,
 the above argument applies if $c_j \sim c_j'$  in $E_{c_j}$ for all $j$,
  even if 
the conditition of being locally-homogeneous is not satisfied.

\begin{sect}{\bf Example.\/} {\rm 
The 
trefoil knot  has
non-trivial invariant with respect to $S_4$ and the cocycle
$\phi={\chi _{0, \,1}} + {\chi _{0, \,T+1}} + {\chi _{1, \,0}} +
{\chi _{ 1, \,T+1}} + {\chi _{T+1, \,0}} + {\chi _{T+1, \,1}}$
over ${\bf Z}_2$
where $S_4={\bf Z}_2[T, T^{-1}]/(T^2 +T +1)$ 
(see \cite{CJKLS}). 
In particular, the color generated by using the
braid form $\sigma_1^2$ with initial colors $0$ and $1,$ gives a
state-sum value of $t.$ 
{}From this braid and the above
construction, we  show that for any $n \in {\bf Z},$ $H_{\rm Q}^2({\bf
Z}_{2n}[T,T^{-1}]/(T^2-T-1), 
{\bf Z}_2) 
\neq 0$ and $H_{\rm Q}^2({\bf
Z}_{2n}[T,T^{-1}]/(T^2+T-1), 
{\bf Z}_2) 
\neq 0$ using the function $\pi:X
\rightarrow S_4, \pi(x)=x \ \  {\mbox{ \rm mod}} \;  2,$  
where $X$ is the quandle for
which we desire the result. For the first case we use virtual
loops colored $2T+2$ and $-4T-1$, and for the latter case we use
loops colored $-2T$ and $-2T+3.$ See Fig.\ref{trefoils}. Finally note
that for the quandle ${\bf Z}_{2n}[T,T^{-1}]/(T^2-T+1),$ the
standard trefoil with (in braid form) initial colors 0,1 and the
cocycle $\pi^\sharp(\phi)$ colors without need for virtual loops and so
$H_{\rm Q}^2({\bf Z}_{2n}[T,T^{-1}]/(T^2-T+1), {\bf Z}_2) \neq 0$.

}\end{sect}

In the spirit of the preceding example, we prove the following.

\begin{sect} {\bf Theorem.\/}
$H_{\rm Q}^2({\bf Z}_2[T, T^{-1}]/(T^2+T+1)^2, {\bf Z}_2) \neq 0$. 
\end{sect}
{\it Proof.\/}
First we describe the virtual knot that we 
use.
Let $\sigma_1$ denote the standard braid generator in $2$-string braid group,
and $v_1$ denote the virtual crossing.
Consider $K_m$ represented by 
$(\sigma_1^2 v_1)^m$.
To compute the colors for $K_m$, 
the
Burau representation is used, 
with the 
matrix 
$B= \left[ \begin{array}{cc} 0&T \\1&1-T \end{array} \right]$ 
replacing $\sigma$ and the  
the permutation matrix 
$P= \left[ \begin{array}{cc} 0&1 \\1&0  \end{array} \right]$ 
replacing $p$. 
Then if the colors assigned to the top two strings on left and right
are $[a,b]$, the color assigned to the strings after the 
sequence $A$ is computed by the matrix multiplication $[a,b]A$.

The matrix corresponding to $K_3$ is 
$$  \left[ \begin{array}{cc} T-T^3 + T^5 -T^6 & T-T^3+T^5\\
1-T+T^3-T^5+T^6 & 1-T+T^3-T^5 \end{array} \right]. $$
Note that $(T^2+T+1)^2$ mod $2$ is $T^4+T^2+1$, and modulo $T^4+T^2+1$
the above matrix is equal to the identity.
Therefore any assignment for the top two strings define a color on 
$K_3$. 

Take for example $[0,1]$ as a color on the top two strings.
The colors assigned to the two strings right above positive 
crossings can be computed as above, and they are 
(starting from the top colors), 
$[0,1]$, $[1,1+T]$, $[0, 1+T]$, $[1+T,T]$,
$[0, T]$, and $[T,1]$, when reduced mod $T^2+T+1$. 
We use the cocycle $\phi'=\pi^\sharp(\phi)$ where
$$ \pi: {\bf Z}_2[T, T^{-1}]/(T^2+T+1)^2 \rightarrow
{\bf Z}_2[T, T^{-1}]/(T^2+T+1)=S_4$$
is the quotient homomorphism and $\phi$ is the cocycle 
described above. 
Hence the state-sum term for this color with the 
cocycle $\phi'$ is $t^3=t$ with $G={\bf Z}_2$ coefficient, 
a nontrivial value. The result follows.
$\Box$

\begin{figure}
\begin{center}
\mbox{
\epsfxsize=2.5in 
\epsfbox{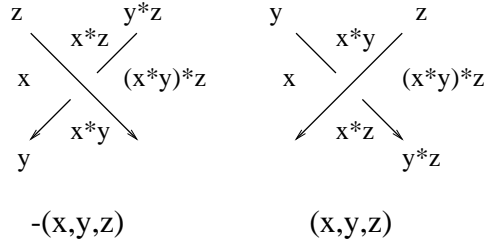}
}
\end{center}
\caption{Shadow colorings of crossings }
\label{shadow}
\end{figure}

In Fig.~\ref{shadow} a local picture for {\it shadow colorings of crossings} is given. 
The regions are colored by quandle elements, as well as over-arcs.
If a region is colored by $x$, an element of a finite quandle $X$, 
then the adjacent region into which the normal of the arc points 
is colored by $x*y$, where $y \in X$ is the color of the arc.
The arcs are colored using the rule defined before.
 Figure~\ref{shadow} shows that 
this rule is well-defined at a crossing. 
Such crossings represent 3-chains as indicated. If a knot or link diagram is 
shadow colored by a  quandle $X$, then the diagram represents a 3-cycle in
$Z^{\rm R}_3(X)$. 
Two isotopic shadow-colored diagrams represent the same homology class.
 We can use shadow colorings
to find non-trivial homology groups as follows.

\begin{figure}
\begin{center}
\mbox{
\epsfxsize=4in 
\epsfbox{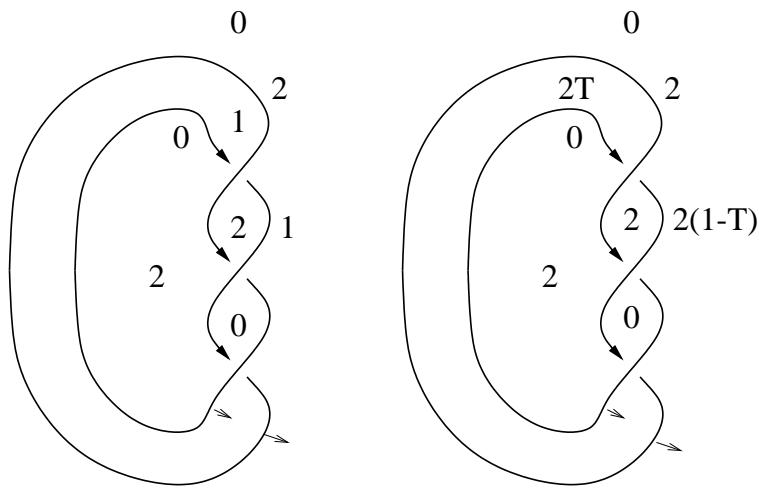}
}
\end{center}
\caption{Shadow colors of trefoil }
\label{treface}
\end{figure}

\begin{sect} {\bf Theorem.\/} 
$H^{\rm R}_3({\bf Z}_3[T, T^{-1}]/(T+1)^2, {\bf Z}_3) \neq 0$. 
\end{sect}
{\it Proof.\/}
We use $\pi: X={\bf Z}_3[T, T^{-1}]/(T+1)^2 \rightarrow R_3$,
$\pi(f(T))=f(-1)$. 
On the left of  Fig.~\ref{treface}, a shadow color by $R_3$ of trefoil is 
depicted, which was used in \cite{RStalk} to show 
that the left and right handed trefoils are distinct.
The diagram on the left of Fig.~\ref{treface} represents the cycle 
$h_0=(2,0,2)+(2,2,1)+(2,1,0)$, 
the class of which is 
a generator 
of $H^{\rm R}_3(R_3, {\bf Z}_3)$.
On the right of Fig.~\ref{treface}, it is shown that the trefoil 
is also colored nontrivially by elements of $X$. 
Let $h_1$ be the class in $Z^{\rm R}_3(X,  {\bf Z}_3)$ represented by 
this face color. Then 
the 3-cycle 
$h_1=(2,0,2)+(2,2,2(1-T))+(2,2(1-T),0)$ 
maps to a non-trivial element in $H_3^{\rm R} (R_3, {\bf Z}_3)$.
Hence $h_1$ is a non-zero element in $H^{\rm R}_3(X, {\bf Z}_3))$.
$\Box$


\end{document}